\documentclass[12pt]{article}
\usepackage{amsmath, amsthm, amssymb}
\usepackage[dvips]{graphics}
\usepackage[usenames,dvipsnames]{xcolor}
\usepackage{subfigure}
\usepackage{multirow}
\usepackage{graphicx} 
\usepackage{bm}     
\usepackage{setspace} 
\usepackage{booktabs} 
\usepackage{verbatim}
\usepackage{appendix} 
\usepackage{natbib} 
\usepackage{changepage}
\usepackage{caption}


\usepackage[bb=boondox,scr=rsfso]{mathalfa} 

\usepackage{soul}

\theoremstyle{plain}

\theoremstyle{definition}

\theoremstyle{remark}

\newcommand{\IG}{\mbox{IG}}
\newcommand{\norm}{\mbox{N}}
\newcommand{\Be}{\mbox{Be}}

\newcommand{\NA}{\mbox{NA}}

\newcommand{\one}{\bm{1}}

\newcommand{\X}{\mathcal{X}}
\newcommand{\GG}{\mathcal{G}}
\newcommand{\CC}{\mathcal{C}}

\renewcommand{\th}{\theta}
\newcommand{\eps}{\epsilon}

\newcommand{\sigs}{\sigma^2}

\newcommand{\xb}{\bm{x}}
\newcommand{\yb}{\bm{y}}

\newcommand{\rs}{r^\star}

\newcommand{\xs}{x^\star}

\newcommand{\xt}{\tilde{x}}

\newcommand{\ind}{\stackrel{\mbox{\scriptsize ind}}{\sim}}

\newcommand{\xmiss}{\xb^{miss}}

\bibpunct{(}{)}{;}{a}{}{,}

\makeatletter
\def\moverlay{\mathpalette\mov@rlay}
\def\mov@rlay#1#2{\leavevmode\vtop{%
   \baselineskip\z@skip \lineskiplimit-\maxdimen
   \ialign{\hfil$\m@th#1##$\hfil\cr#2\crcr}}}
\newcommand{\charfusion}[3][\mathord]{
    #1{\ifx#1\mathop\vphantom{#2}\fi
        \mathpalette\mov@rlay{#2\cr#3}
      }
    \ifx#1\mathop\expandafter\displaylimits\fi}
\makeatother

\begin{document}
\def\spacingset#1{\renewcommand{\baselinestretch}%
{#1}\small\normalsize} \spacingset{1}

  \title{\bf Regression with Variable Dimension Covariates}
  \author{Peter Mueller \\ Department of Statistics \& Data Science \\ University of Texas at Austin \\
    and \\
 Fernando Andr\'es Quintana \\ Departamento de Estad\'{i}stica, \\ Pontificia Universidad Cat\'{o}lica de Chile, Santiago
 \\ and Millennium Nucleus Center for the Discovery of Structures in Complex Data \\
  	and \\ Garritt L. Page \\
    Department of Statistics\\ Brigham Young University, Provo, Utah
    }
  \maketitle

\begin{abstract}
	 Regression is one of the most fundamental statistical inference
	problems.  A broad definition of regression problems is as estimation
	of the distribution of an outcome using a family of probability
	models indexed by covariates.
	Despite the ubiquitous nature of regression problems and the
	abundance of related methods and results there is a surprising gap
	in the literature.
	There are no well established methods for regression with a varying
	dimension covariate vectors, despite the common occurrence of such
	problems.
	In this paper we review some recent related papers proposing 
	varying dimension regression by way of random	
	partitions.
\end{abstract}


\noindent%
{\it Keywords: }{density regression, clustering, partition, missing data
\vfill

\newpage
\spacingset{1.5} 

\section{Introduction}\label{intro}
We discuss approaches for Bayesian inference for regression with
varying dimension covariate vectors.
We review a sequence of recent papers that develop an approach based
on random partitions and a cluster-specific outcome model.
The random partition of experimental units is set up in a way that allows the use of
any available subset of a list of covariates. This formalizes the
intuitive notion of clustering experimental units on the basis of
available information, as it is commonly practiced in everyday
problems.  The resulting scheme is a nonparametric Bayesian regression
that works with available covariates for each experimental unit,
allowing any subset of a full covariate vector.
The only major assumption is that covariates are missing at random. No
further structural assumptions are needed.

Consider the generic regression problem of explaining an outcome
$y_i$ as a function of a covariate $\xb_i \in \X$.
For the moment we assume $y_i \in \Re$ and $\xb_i \in \Re^p$, 
and state the regression problem as
$
y_i = f(\xb_i) + \epsilon_i
$,
$i=1,\ldots,n$.
Here $f$ is an unknown centering function and $\eps_i$ are residuals,
usually assumed to be independent.
In traditional parametric regression the
function $f$ and the residual distribution are indexed by a finite
dimensional parameter vector $\th$.
Without the restriction to finite dimensional $\th$ we are led to
nonparametric extensions of regression problem.
In the most general case, without finite parametric model for either,
the problem is characterized as
\begin{equation}
	y_i \mid \GG, \xb_i \ind G_{\xb_i}
\end{equation}  
with a prior model
$\pi(\GG)$ on the family 
$\GG=\{G_{\xb},\; \xb \in \X\}$
of outcome distributions indexed by $x$.
Prior probability models for random distributions (and families of
random distributions) are known as nonparametric Bayesian (BNP) models. See
\cite{ghosal2017fundamentals} for an extensive discussion of
underlying models and theory. 
The most widely used BNP model remains the 
Dirichlet process \citep{ferguson1973bayesian} and its variations and
extensions. An early careful discussion of the Dirichlet process 
and its properties appears in
\cite{basu2011note}. \cite{sethuraman2011commentary} provides 
delightful comments on the background and history of that contribution.

Implicit in the previous description of regression is the assumption of complete covariate
vectors $\xb_i$.
Most discussions of regression, including BNP regression, follow this
assumption.
A generic solution strategy, of course, is to treat incomplete
covariate vectors as a missing data problem and impute the missing
values.
Many model-based methods have been developed to make better use of
information with missingness, including maximum likelihood (ML)
methods, multiple imputation (MI) methods, weighted estimating
equation (WEE) methods, and fully Bayesian (FB) methods.  Detailed
reviews of these methods appear in
\cite{little1992regression},
\cite{horton1999maximum},
\cite{schafer2002missing}, and
\cite{ibrahim1999missing}.
The ML, MI, WEE and FB approaches
require an exposure model $p(x_i \mid \alpha)$ for covariates
in addition to an outcome model
$p(y_i \mid x_i, \th)$.
Here $\th$ indexes the outcome model and $\alpha$ denotes the set of
parameters for the exposure model.
For example,
\cite{lipsitz1996conditional} and
\cite{ibrahim1999missing}
construct $p(x_i \mid \alpha)$ as
a product of one-dimensional conditional distributions:
$$p(x_{i1} \mid \alpha_1), \ldots,
p(x_{i,p-1} \mid x_{i1}, \ldots, x_{i,p-2}, \alpha_{p-1}) 
\cdot
p(x_{i,p} \mid x_{i1}, \ldots, x_{i,p-1}, \alpha_{p}).
$$
Specifying a probability model for $x_i$
is intuitively appealing and usually convenient to implement, but
becomes challenging for high- to moderate dimensional covariates.
Some approaches address this challenge using
simultaneous variable selection or tree-based methods.
For example \cite{jiang2022adaptive} use iteratively updated missing
data and hyper-parameters. Specifically, they consider a combination
of $L^1$ regularization with variable selection methods and a covariate
imputation scheme based on a stochastic approximation to the
expectation-maximization algorithm (SAEM). Alternatively, \cite{mercaldo&bluem:2018} consider a strategy based on
pattern submodels, that is, a set of submodels for every missing data pattern and which are fit using data only from that particular pattern.

While these are
valid and principled approaches, and very natural in the
context of simulation based Bayesian inference, it could be argued that
in everyday regression and decision problems agents proceed in a more
parsimonious manner. For example, a clinician considering treatment
options would consider possible outcomes based on all available
patient covariates, using available information, but not imputing
missing information (unless some evidence gives rise to suspect
informative missingness, like lab values below detection limits etc.).
Grouping people in a social context we routinely use available
information, grouping, for example, speakers at a conference with
respect to some characteristics, and quite possibly missing many variables that
could be helpful in clustering speakers if we knew.
In this paper we review a recently introduced approach to formalize
this process as nonparametric Bayesian regression based on random
partitions.

All reviewed approaches are based on random partition models (for
units and for missingness patterns). That is, probability models for
cluster arrangements. We build on the product partition model
\citep{hartigan1990partition,
	quintana&loschi&page:2018}.
The PPM model has been used for BNP data analysis
before in many other contexts, including estimation of normal means
\citep{CRO97}, identification of changepoints in time series
\citep{loschi2005multiple}, and disease mapping \citep{hegarty&barry:08}.
In particular, the popular Chinese restaurant process, a term
introduced by Jim Pitman and Lester Dubins
\citep[see, e.g.][]{aldous:85,PI96}
fits into the PPM framework too
\citep{quintana2003bayesian}.

In Section \ref{sec:PPMx} we introduce the basic model based on a
random partition of experimental units. 
Section \ref{sec:synthetic} discusses an application to creating
synthetic matching (patient) populations in the presence of variable
dimension covariate vectors.
In Section \ref{sec:PPMx-regr} we extend the basic model by
introducing cluster-specific regression sub-models.

\section{Regression with variable dimension covariates
	using random partitions}
\label{sec:PPMx}
In \cite{page2022clustering} we introduce an approach using regression
based on a random partition.
In words, we introduce a partition $\CC=\{C_1,\ldots,C_K\}$ of
experimental units $[n]=\{1,\ldots,n\}$ based on covariates $\xb_i$
and cluster-specific parameters $\th_k$ for an outcome model. Here the
prior on $\CC$ is such that units $i,i'$ with more similar covariates
are more likely to co-cluster. This is achieved using the PPMx model
introduced in \cite{PPMxMullerQuintanaRosner}. The latter is a prior
model $p(\CC \mid \xb)$ that is constructed to favor clusters with
similar covariates. The desired random partition $p(\CC \mid \xb)$ is
defined as a product partition model 
(PPM) \citep{hartigan1990partition} using a cohesion function with an
additional factor that is designed to favor similar covariates for all
units in the cluster. In this judgment of similarity, for each
covariate only the units that report values
are employed. Missing covariate values are simply skipped. 

We introduce some notation for a formal description.
Let $r_{ij} \in \{0,1\}$ denote an indicator for variable $j$ for
subject $i$ being reported. That is, $1-r_{ij}$ is an indicator for
missingness, 
let $O_{kj} = \{i:\; i \in C_k \mbox{ and } r_{ij} =1\}$ be the set
of all units in $C_k$ with available data for the $j$-th covariate,
and let 
$\xs_{kj} = (x_{ij};\; i \in O_{kj})$ denote the reported covariates
$x_{ij}$ grouped by cluster, and $\xs_k=\{\xs_{k1},\ldots,\xs_{kp}\}$.
The PPMx model defines $p(\CC \mid x)$ for a covariate-dependent
random partition as 
\begin{equation}
	p(\CC \mid \xb) \propto \prod_{k=1}^K c(C_k)\, g(\xs_k)
	\mbox{ with }
	g(\xs_k) = \prod_j g_j(\xs_{kj})
	\label{PPMx}
\end{equation}
where $g_j(\xs_{kj})$ is a function (``similarity function'')
that scores the similarity of the values in $\xs_{kj}$, similar to a
purity function in hierarchical clustering
\citep{InfRetrieval}.
It returns maximum values for all equal $x_{ij}$, $i \in C_k$,
and low values for very diverse values.
A convenient formalization is as a marginal probability in a conjugate
model, as
\begin{equation} \label{g*}
g(\xs_{kj}) =
\int \prod_{i \in O_{kj}} q(x_{ij} \mid \xi_{kj})\, dq(\xi_{kj}).
\end{equation}
Model $q(\cdot)$ in \eqref{g*} is said to be auxiliary in the sense
that it is only used for computational convenience,
without any notion of modeling $x$.
Let $N(x;\, m,V)$ denote a normal p.d.f. with moments $(m,V)$,
evaluated at $x$, and similarly, let $\IG(x;\; a,b)$ denote an inverse
gamma pdf (with mean $b/(a-1)$) evaluated for $x$.
For example, for a continuous variables $x_{ij}$ one could use
\begin{eqnarray*}
q(x_{ij} \mid \xi_{kj}=(\mu_{kj}, \sigs_{kj})) & = &N(x_{ij};\; \mu_{kj},
\sigs_{kj})\\
q(\xi_{kj}) & = & \IG(\sigs_{kj};\; a_{kj},b_{kj})\, N(\mu_{kj};\; 0,
c_{kj}\sigs_{kj})
\end{eqnarray*}
Here $(a_{kj},b_{kj},c_{kj})$ are fixed hyperparameters, chosen to
reflect the range of plausible values for the $j-$th covariate
and the desired characterization of similarity.
The definition of $g(\cdot)$ as a marginal distribution under the
auxiliary model $q$ exploits the fact that the marginal -- in this case a
version of a multivariate t-distribution -- is most peaked for very
similar $x_{ij}$.
Similarly, for binary variables we use the marginal beta-binomial
distribution of the binary outcomes.
However, there is no notion of modeling a covariate distribution. The
use of the auxiliary model $q$ is merely for easy calculus to evaluate
$g(\cdot)$.
Any alternative function could be used. For example, for catagorical
data one could use the relative frequency of the most common value in
each cluster.

The random partition model \eqref{PPMx} is then completed with an
outcome model
\begin{equation}
p(y_i \mid i \in C_k, \th) \sim p(y_i \mid \th_k).
\label{pyth}
\end{equation}
Let $D=(x_i,y_i;\; i=1,\ldots,n)$ denote the observed data.
Models
\eqref{PPMx} and \eqref{pyth} together imply a predictive distribution $p(y_{n+1} \mid
x_{n+1}=x, D)$ for a future outcome as a function of covariates as
\begin{equation}
p(y_{n+1} \mid x, D) =
\int
\sum_{k=1}^K p(y_{n+1} \mid \th_k)\, p(n+1 \in C_k \mid x, \CC, D)\;
dp(\CC, \th \mid D).
\label{PPMxy}
\end{equation}
Here $p(\CC,\th \mid D)$ refers to the posterior probability model for
the random partition and the cluster-specific outcome parameters
$\th_k$, and $p(n+1 \in C_k \mid x,\CC,D)$ is the 
probabilty of adding a new, $(n+1)-$st unit with $x_{n+1}=x$
to cluster $C_k$.
Defining similarity functions with an
auxiliary probability model as in \eqref{g*}
has the appealing property of rendering
a sample size consistent model for $\CC$, i.e., the model for $n$ units arises from
that for $n+1$ by marginalizing the last one.
See the discussion in \cite{PPMxMullerQuintanaRosner}.
%
In words, the prediction for a future outcome is obtained by first
allocating the new unit in one of the (imputed) clusters $C_k$,
favoring clusters with similar covariates; given the cluster
membership the prediction is then based on the cluster-specific outcome
parameter $\th_k$. 
The reported regression $p(y_{n+1} \mid x, D)$ averages w.r.t. the
posterior on $\CC$ and $\th$.

An important feature of $p(y_{n+1} \mid x, D)$ is that it is
well-defined for any subset of available covariates in
$x=(x_1,\ldots,x_p)$.
This is because $p(n+1 \in C_k \mid x_{n+1},\CC,D)$ uses the available
covariates only.
From \eqref{PPMx} and \eqref{g*} we have
\begin{multline}
p(n+1 \in C_k \mid x_{n+1},\CC,D)
\propto
\frac{c(C_k \cup \{n+1\})}{c(C_k)} \times \\
\prod_{j:\; r_{n+1,j}=1}
\frac
{\int \prod_{i \in O_{kj}\cup \{n+1\}} q(x_{ij} \mid \xi_{kj})\, dq(\xi_{kj})}
{\int \prod_{i \in O_{kj}\phantom{\cup \{n+1\}}} q(x_{ij} \mid
	\xi_{kj})\, dq(\xi_{kj})}\; 
\label{eq:pred}
\end{multline}
This is illustrated in Figure \ref{fig:pred} by showing the predictive
$p(y_{n+1} \mid x_{n+1}, \CC, \th)$ for $p=2$ covariates,
that is \eqref{PPMxy} before posterior averaging w.r.t. $\CC,\th$.
The figure shows the regression for
complete data $x_{n+1}=(x_{n+1,1},x_{n+1,2})$ (black surface),
for one missing covariate, $x_{n+1}=(x_{n+1,1}, \NA)$ (red curve in
the xz-plane), and for all missing covariates, $x_{n+1}=\NA$ (green
bullet on the z-axis). 
\begin{figure}[hbt]
\includegraphics[width=.8\textwidth]{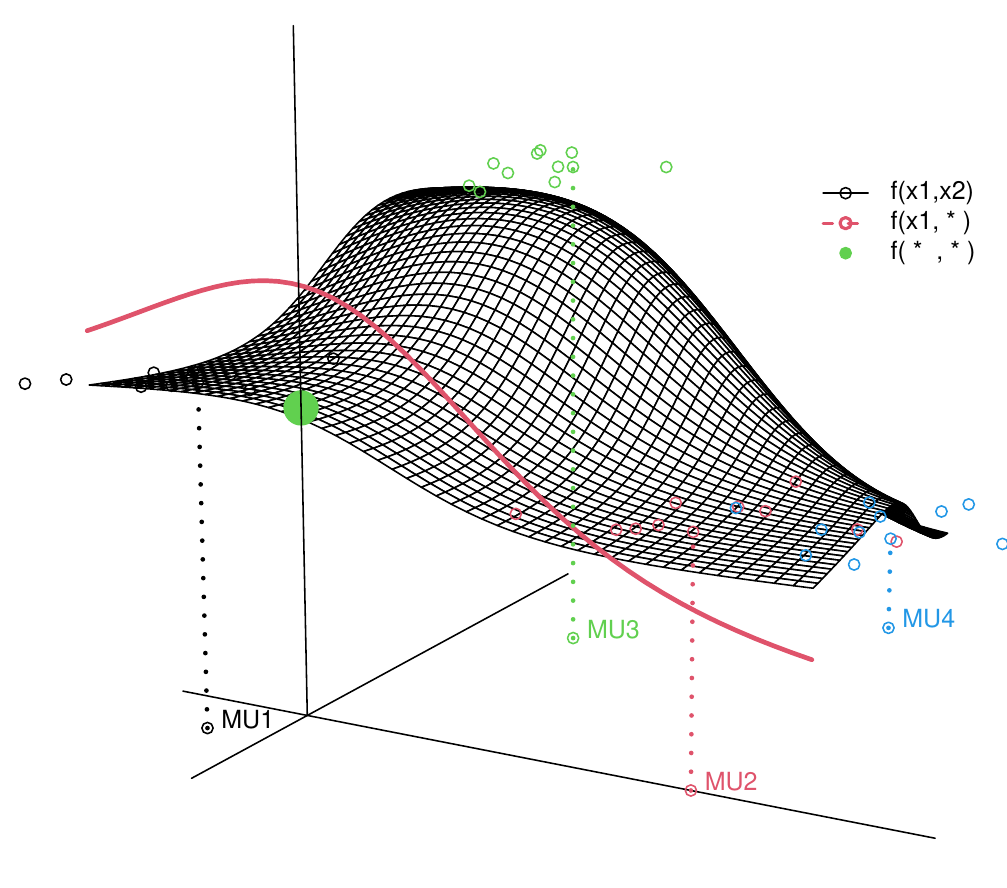} 
\caption{The figure shows the regression $f(y_{n+1} \mid x, D)$ for a
complete covariate vector $x=(x_1,x_2)$ as the black surface; for
$x=(x_1,N/A)$ as the red curve (in the xz-plane), and for all
missing covariates as the green bullet (on the z-axis).
The data are shown as small circles, with a random partition into
four clusters $C_1, \ldots, C_4$, shown in black, red, green and
blue. The cluster centers are indicated as $\mu_1, \ldots, \mu_4$.}
\label{fig:pred}
\end{figure}

\paragraph*{\bf\em Example 1: Simulation.} 
We generate (complete) data $(y_i,\xt_i)$, $i=1,\ldots,n=160$ for $p=4$
covariates, $K=8$ distinct missingness patterns, and $m=20$
observations per pattern, using $y_i \sim N(\xt_i'\beta,\sigs)$, with
$\xt_i = (\xt_{ij};\; j=0,\ldots,p)$, $\beta=(2, 1.4, 1, 0.1, 2)$,
with an intercept
$\xt_{i0}=1$,
a rescaled beta distribution
$\xt_{i1} \sim 0.5 \cdot \Be(4,1)$,
a correlated second covariate
$\xt_{i2} \sim \xt_{i1} + N(\xt_{i1},1)$,
a mixture of two normals
$\xt_{i3} \sim 0.3\cdot \norm(-3,1) + 0.7 \cdot \norm(3,1)$,
and a rescaled bimodal beta 
$\xt_{i4} \sim 5 \cdot \Be(0.3,0.3)$.
We implement posterior inference for $D=\{(y_i,x_i)\}$
with incomplete covariate vectors $x_{ij} = r_{ij}\xt_{ij}$, using
$r_i = \one$ for the first $m=20$ observations, and
$r_i = \rs_k$ for 20 observations each, for $k=2,\ldots,8$.
Here
$\rs_2=(0,1,1,1)$, $\rs_3=(1,0,1,1)$, $\rs_4=(1,1,0,1)$,
$\rs_5=(1,1,1,0)$, $\rs_6=(0,0,1,1)$, $\rs_7=(0,1,0,1)$ and $\rs_8=(1,0,1,0)$.

The described data generation scheme was used to generate 100
datasets of $n=160$ observations each.
For each data set 10\% of the observations (16) were
randomly selected to comprise the testing dataset while the remaining
144 comprised the training data set.  
We then carried out regression for each data set using the following approaches:
(1) VDReg, as described in Section \ref{sec:PPMx}. More specifically,
we used model (10) in \cite{page2022clustering}. Inference under the
VDReg model is implemented using the {\tt
ppmSuite} package (\citealt{ppmSuite}) in {\tt R}; 
(2) BART for regression with missing covariates as 
introduced in \cite{kapelner&bleich:2015}. The method
uses missing covariates to inform splitting decisions when growing
regression trees and is fit using the {\tt bartMachine} package; 
(3) PSM, the pattern submodel approach proposed in
\cite{mercaldo&bluem:2018}. The approach specifies a separate regression model
(all of the same class) for each missing pattern and is fit using
sotfware available at {\tt
https://github.com/sarahmercaldo/MissingDataAndPrediction}.

For each generated dataset, we computed the mean squared prediction
error (MSPE) based on the 16 out-of-sample predictions for the test
set.  
We then recorded the average for each approach across the 100 synthetic
datasets.  
These averages are shown in Table \ref{tab:mspe}.  
The VDReg procedure and BART have similar out-of-sample
prediction rates, while the PSM reported the largest MSPE amongst the three
considered procedures.

\begin{table} 
\caption{MSPE values averaged over the 100 simulated datasets for each procedure}
\label{tab:mspe}
\begin{center}
\begin{tabular}{l  ccc} \toprule
	Procedure & VDReg  & BART & PSM\\ \midrule
	MSPE & 6.52  &  6.23  & 7.11\\ \bottomrule
\end{tabular}
\end{center}
\end{table}

\section{Synthetic matching populations with
variable dimension data vectors}
\label{sec:synthetic}
In \cite{chandra2021bayesian} we use the approach of Section
\ref{sec:PPMx} to generate a synthetic control cohort for a single-arm
treatment only clinical study.
Let then $D_1 =\{(y_{1i}, \xb_{1i});\; i=1,\ldots,n_1\}$ denote the data in a
single-arm treatment only study with $n_1$ enrolled patients, with
baseline covariates $\xb_{1i}$ and outcomes $y_{1i}$. In the
motivating application the study is a clinical trial for glioblastoma (GBM)
patients, the outcomes $y_{1i}$ are overall survival (with censoring),
and the covariates $\xb_{1i}$ are $p=10$ important baseline covariates
that are commonly used in GBM studies.
While single-arm trials are common in early phase GBM studies, the
desireable gold standard for clinical studies is still a randomized
clinical trial with random assignment of patients to treatment and
control arms. One of the reasons for using treatment-only trials in
GBM are difficulties in patient recruitment, and the lack of an
effective active control.
The intention is then to use available historical data from earlier
studies to construct a synthetic control cohort.
Let  $D_2 = \{(y_{2i}, \xb_{2i});\; i=2,\ldots,n_2\}$ denote the
historical data.
A critical feature of this approach is
that historical patients should be selected such that the two patient
cohorts can be considered to be equivalent, i.e., matching
distributions of baseline covariates.

\cite{chandra2021bayesian} propose an approach fitting the PPMx model
from Section \ref{sec:PPMx} to $D_1$. Let $\CC_1 =
\{C_{11},\ldots,C_{1K}\}$ denote a random partition of $[n_1]$, and
let $\th_{1k}$ denote parameters for the cluster-specific outcome model.
Data $D_2$ is then partitioned to create clusters
$C_{21}, \ldots, C_{2K}$ matching $\CC_1$, plus additional clusters if
needed, and introducing $\th_{2k}$ as cluster-specific parameters for
an outcome model in $D_2$.
Assuming for the moment that $n_2$ is usually much larger than $n_1$,
we can constrain the model to $|C_{2k}| \ge |C_{1k}|$. Dropping patients in the
additional clusters and thinning out $C_{2k}$, $k=1,\ldots,K$ to match
the cluster sizes $|C_{1k}|$ one can then achieve matching
populations. The actual implementation works with weights instead of
dropping data points.

There are two important features in this process. First, historical
data usually includes a good number of missing data. Considering the
carefully controlled context of clinical studies one can assume
missing at random (for example, some variables were not recorded in an
earlier study). The described approach implements inference without
the need to impute such missing data.
Second, the model implements BNP regression (with variable subsets of
covariates), with a pair of outcome model parameters
$(\th_{1k},\th_{2k})$ in each cluster. This allows one to define
cluster-specific treatment effects $\delta_k=d(\th_{1k},\th_{2k})$,
using an appropriate function $d$. For example, if $\th_{sk}$ has the
interpretation of a mean-outcome, one could use
$d(\th_1,\th_2)=(\th_1-\th_2)$.
Cluster-specific $\delta_k$ can be averaged to define an
overall treatment effect $\Delta$, with a full probabilistic description of
uncertainties. In particular, reported inference on $\Delta$ averages
over the random partition and all unknown parameters.

\section{Including cluster-specific regression}
\label{sec:PPMx-regr}

Motivated by the goal of improving predictive capabilities,
\cite{heiner&page&quintana} generalized the approach from Section
\ref{sec:PPMx} by allowing the cluster-specific outcome model
$p(y_i\mid c_i=k, \th_k)$ to be now specified as a regression
$p(y_i \mid \xb_i, c_i=k, \th_k)$.
That is, a regression model $p(y_i \mid \xb_i, c_i=k, \th_k)$ replaces
the outcome model \eqref{pyth} with a local regression.
Following similar considerations, \cite{friedberg2020local} find it
useful and advantageous to incorporate local predictors in the
context of random forest models.
However, in the context of missing covariates, 
this approach faces the practical problem of requiring all
covariates in the local regression model, including missing ones.
In \cite{heiner&page&quintana} this problem was addressed by noting that
analytically integrating
out the missing values in $\xb$ w.r.t. the auxiliary model in
\eqref{g*} yields 
the same distribution of $(\yb, \rho \mid \xb)$ as would
arise from modeling $(\yb \mid \CC)$ with $(\rho \mid
\xb)$ using $g$ in \eqref{g*}. In other words, skipping
over missing covariates from the similarity scores is, under certain
conditions equivalent to integrating them out of a PPM that treats
$\xb$ as random. 
We refer to this step as ``projection''
(understood here as a synonym of ``marginalization''), indicating that $\xb$ is
not modeled in any way.  Importantly, the scheme
still entirely avoids imputations.

We can carry this idea a bit further and relax the independence of
$\yb$ and $\xb$, still obtaining an analytically tractable
scenario. To explain the idea, we momentarily drop the subject index
$i$. Let $q_j\left(x_j\right)=N\left(x_j ; \mu_j^{(x)}, \sigma_j^{(x)
2}\right)$ denote the auxiliary model for covariate $j=1, \ldots, p$,
and $y \mid \xb \sim N\left(\mu+\sum_{j=1}^p \beta_j z_j,
\sigma^2\right)$, where $z_j=\left(x_j-\mu_j^{(x)}\right) /
\sigma_j^{(x)}$. Integrating the joint density with respect to
the missing values $\xmiss$ in $\xb$ yields
\begin{equation}
\int p(y \mid \xb) \prod_{j}^p q_j\left(x_j\right) \mathrm{d} \xmiss=
N(m, V)\, \prod_{j:\,r_j=1} q_j\left(x_j\right), 
\label{eq:condmodel}
\end{equation}
with 
$m = \sum_{j:\,r_j=1} \beta_j z_j$ and
$V= \sigma^2+\sum_{j:\, r_j=0} \beta_j^2$.
The introduction of the centered and scaled covariates $z_j$ 
stabilizes the mean and simplifies the expression for the
inflated variance of the conditional distribution of $y$.

Note that \eqref{eq:condmodel} can be stated without reference to the missing
values in $\xb$. Also, the $\left\{(\mu_j^{(x)},
\sigma_j^{(x)})\right\}$ parameters play no role in the actual
model, and can be replaced by conveniently chosen 
plug-in alternatives, such as posterior means and variances under
customary conjugate alternatives, say $\hat{\mu}_j^{(x)}$ and
$\hat{\sigma}_j^{(x) 2}$.

Putting all of this together, the variable dimension covariate model with
local linear regression (VDLReg) poses a likelihood specification as follows
for $i=1, \ldots, n$, $j=1, \ldots, p$, and $k=1, \ldots, K$:
\begin{multline}\label{eq:VDLReg}
y_i \mid c_i=k, \theta_k\ind
N(m,V)\\
m=\mu_k+\sum_{j:\, r_{ij}=1} \beta_{kj} z_{ij},~~
V= \sigma^2_k+\sum_{j:\, r_{ij}=0} \beta_{kj}^2
\end{multline}
where $\theta_k=(\mu_k,\beta_{k1},\ldots,\beta_{kp},\sigma_k^2)$.

One additional aspect of model \eqref{eq:VDLReg} is the increased variance
that comes from projecting the missing values. This
could limit predictive performance.
In \cite{heiner&page&quintana} we addressed this problem
by aggressively shrinking the regression coefficients
$\left\{(\beta_{k1}, \ldots, \beta_{kp})\right\}$ with the adoption of
a Dirichlet-Laplace prior \citep{bhattacharya2015dirlap} at the cluster level.

We added the VDLReg approach to the simluaiton study summarized in Table
\ref{tab:mspe}.
We implemented inference under the  VDLReg approach for the same 100
datasets that were generated in Example 1 in Section \ref{sec:PPMx}.  In addition to
\ref{eq:VDLReg},  details associated with the model that was fit are
provided in equation  (4) of \cite{heiner&page&quintana}.  We fit the
VDLReg  procedure using {\tt Julia} code available at {\tt
https://github.com/mheiner/ProductPartitionModels.jl}.  The MSPE
based on the same testing observations turned out to be 5.81 for
VDLReg, that is, the smallest among all four considered approaches.

\section{Conclusion}
We reviewed some approaches to implement regression and prediction
with varying dimension covariate vectors, as it is commonly
done in everyday decision making, but curiously overlooked in the
statistics literature.
The proposed approaches are based on regression by clustering. That
is, we first partition experimental units into subgroups that are
judged similar on the basis of available covariates, and then assume
an outcome model for each cluster. The important detail here is that
the random clustering is set up on the basis of all available
covariates, without imputing missing covariates.
This brief description also already points to the main
limitation. Informative missingness makes the approach invalid.

Also the construction of a suitable similarity function is potentially
challenging. Using the default computation-friendly solution as the
marginal under a conjugate auxiliary model is convenient, but leaves
inference actually identical to what it would be under an extended
outcome of response and covariates combined (as discussed in Section
4).
But the framework is more general, and allows for any desired
similarity function, at the cost of less computation-efficient
posterior simulation.
However, if an application suggests problem-specific similarity functions the
additional computational effort is a reasonable cost for 
being able to accommodate relevant expert judgment and decision maker
preferences.

Finally, as briefly mentioned before, the use of local
cluster-specific regression models in VDLReg highlights the similarity
with tree-based regression, which might use local regression in each
leave of the tree. The main difference is that tree-based methods work
with partitions of the covariate space, usually using rectangular
subsets defined by sequences of thresholds. In contrast, the approach in VDLReg
allows for more general random partitions.

In summary, the discussed approaches are most suitable for problems
with massive missing data, with missingness for well-understood
reasons and non-informative, and informed expert judgment on relevant
similarity of experimental units.
The non-parametric BNP nature of the approach is attractive when
biases due to parametric assumptions are problematic, as BNP models
are usually ``always right'' (in the formal sense of full prior
support). This makes VD(L)Reg particularly useful for applications in
biomedical problems. We discussed a typical application in Section 3,
and believe the approach could be useful in many more problems related
to the design and data analysis for clinical studies.
\label{sec:concl}

\section*{Acknowledgments}

Fernando Quintana was partially supported by grant FONDECYT 1220017,
Peter M\"uller was partially supported by NSF under grant NSF/DMS 1952679. 
\clearpage
\bibliographystyle{dcu}
\bibliography{reference}
\end{document}